\documentclass{elsart}

\usepackage{amssymb,amsmath}
\usepackage{url}

\newcommand{\bbOmega}{\boldsymbol \Omega}
\newcommand{\Proj}{\mathbb{P}}
\newcommand{\Fp}{\mathbb{F}_p}
\newcommand{\Fq}{\mathbb{F}_q}
\newcommand{\Fqr}{\mathbb{F}_{q^r}}
\newcommand{\Fqrho}{\mathbb{F}_{q^\rho}}
\newcommand{\parameter}{\psi}
\newcommand{\parametercurve}{\lambda}
\newcommand{\parametercurvebis}{\mu}
\newcommand{\expcurvea}{a}
\newcommand{\expcurveb}{b}
\newcommand{\quintic}[1]{\mathcal{M}_{#1}}
\newcommand{\mirrorquintic}[1]{\mathcal{W}_{#1}}
\newcommand{\curveA}[1]{\mathcal{A}_{#1}}
\newcommand{\curveB}[1]{\mathcal{B}_{#1}}
\newcommand{\curveC}[1][]{\mathcal{C}_{#1}}
\newcommand{\zetafunction}[2]{Z_{#1}\!\left(#2\right)}
\newcommand{\card}[1]{|#1|}
\newcommand{\set}[2]{\{#1~|~#2\}}
\newcommand{\Fqhat}{\widehat{\mathbb{F}_q^*}}
\newcommand{\nbzeros}[1]{N_{#1}}

\newcommand{\caradd}{\varphi}
\newcommand{\caraddbis}{\varphi '}
\newcommand{\carmult}{\chi}
\newcommand{\carmultbis}{\eta}
\newcommand{\carmulttriv}{\boldsymbol{1}}
\newcommand{\carmultquad}{\varepsilon}
\newcommand{\gausssum}[1]{G(\caradd,#1)}
\newcommand{\gausssumadd}[2]{G(#1,#2)}
\newcommand{\jacobisum}{J}

\DeclareMathOperator{\Tr}{Tr}
\newcommand{\TrFqFp}{\Tr_{\Fq/\Fp}}
\newcommand{\TrFqrFq}{\Tr_{\Fqr/\Fq}}
\newcommand{\NFqrFq}{N_{\Fqr/\Fq}}

\begin{document}

\begin{frontmatter}



\title{On the zeta function of a family of quintics}


\author{Philippe Goutet}

\address{IMJ, case 247, 4 place Jussieu, 75252 Paris Cedex 05, France}

\ead{goutet@math.jussieu.fr}

\begin{abstract}
In this article, we give a proof of the link between the zeta function of two families of hypergeometric curves and the zeta function of a family of quintics that was observed numerically by Candelas, de la Ossa, and Rodriguez Villegas. The method we use is based on formulas of Koblitz and various Gauss sums identities; it does not give any geometric information on the link.
\end{abstract}

\end{frontmatter}

\section{Introduction}
\label{intro}

Let $\Fq$ be a finite field of characteristic $p \neq 5$. In all this article, $\parameter$ will be an element of $\Fq$ and $\quintic{\parameter}$ the hypersurface of $\Proj^4_{\Fq}$ defined by
\[x_1^5 + \dots + x_5^5 - 5 \parameter x_1 \dots x_5 = 0.\]
It is non-singular if and only if $\parameter^5 \neq 1$. We denote by $\card{\quintic{\parameter}(\Fqr)}$ the number of points of $\quintic{\parameter}$ over an extension $\Fqr$ of degree $r$ of $\Fq$. The zeta function of $\quintic{\parameter}$ is
\[\zetafunction{\quintic{\parameter}/\Fq}{t} = \exp\left(\sum_{r=1}^{+\infty}{\card{\quintic{\parameter}(\Fqr)} \frac{t^r}{r}}\right).\]
Candelas, de la Ossa, and Rodriguez Villegas have given in \cite{CdlORV.I,CdlORV.II} an explicit formula for $\card{\quintic{\parameter}(\Fqr)}$ in terms of Gauss sums. The formula takes the form
\[\begin{split}
\card{\quintic{\parameter}(\Fqr)} = 1 & + q^r + q^{2r} + q^{3r} \\ & - N_w(q^r) - 10q^r N_a(q^r) - 15q^r N_b(q^r) + 24 N_\textup{sing}(q^r).
\end{split}\]
When $\quintic{\parameter}$ is non singular (that is, when $\parameter^5 \neq 1$), the term $N_\textup{sing}(q^r)$ is zero (see lemma~\ref{result:nb.zeros.sing} page~\pageref{result:nb.zeros.sing} below). In this case, we have
\begin{equation}\label{formula:zeta.quintic}
\zetafunction{\quintic{\parameter}/\Fq}{t} = \frac{P_w(t) P_a(qt)^{10} P_b(qt)^{15}}{(1-t)(1-qt)(1-q^2t)(1-q^3t)},
\end{equation}
where $P_w$ is the formal series $\exp\left(\sum_{r=1}^{+\infty}{N_w(q^r) \frac{t^r}{r}}\right)$, $P_a$ and $P_b$ being defined in a similar way. It is possible to associate to $\quintic{\parameter}$ a ``mirror variety'' $\mirrorquintic{\parameter}$ (obtained by quotient and then resolution of singularities, see \cite[\S 10 p.~124]{CdlORV.II}). When $\parameter \neq 0$, Candelas, de la Ossa, and Rodriguez Villegas \cite[\S 14 p.~149]{CdlORV.II} show that
\[\zetafunction{\mirrorquintic{\parameter}/\Fq}{t} = \frac{P_w(t)}{(1-t)(1-qt)^{101}(1-q^2t)^{101}(1-q^3t)}.\]
As the Betti numbers of $\mirrorquintic{\parameter}$ are (1,0,101,4,101,0,1), we know from the Weil conjectures that $P_w$ is a polynomial of degree 4. When $\parameter = 0$, the above formula for $\zetafunction{\mirrorquintic{\parameter}/\Fq}{t}$ is still valid, but it is also possible to give a hypergeometric interpretation of $P_w(t)$; more precisely, we will show (see remark~\ref{remark:H0} page~\pageref{remark:H0}) that $P_w(t)$ is the numerator of the zeta function of the hypergeometric hypersurface
\[H_0 : y^5 = x_1 x_2 x_3(1-x_1-x_2-x_3).\]

Concerning $P_a$ and $P_b$, define two affine curves by $\curveA{\parameter}$~:~$y^5 = x^2(1-x)^3(x-\parameter^5)^2$ and $\curveB{\parameter}$~:~$y^5 = x^2(1-x)^4(x-\parameter^5)$. We write the corresponding zeta functions as
\[\zetafunction{\curveA{\parameter}/\Fq}{t} = \frac{P_{\curveA{\parameter}/\Fq}(t)}{1-qt} \quad \text{ and } \quad \zetafunction{\curveB{\parameter}/\Fq}{t} = \frac{P_{\curveB{\parameter}/\Fq}(t)}{1-qt},\]
where $P_{\curveA{\parameter}/\Fq}$ and $P_{\curveB{\parameter}/\Fq}$ are of degree $8$ if $\parameter \neq 0$ and of degree $4$ if $\parameter = 0$. On the basis of numerical observations (made in the case $q=p$ for $p \leq 101$), Candelas, de la Ossa, and Rodriguez Villegas conjecture that $P_a = P_{\curveA{\parameter}/\Fq}$ and $P_b = P_{\curveB{\parameter}/\Fq}$ if $\parameter \neq 0$ and that $P_a = P_{\curveA{\parameter}/\Fq}^2$ and $P_b = P_{\curveB{\parameter}/\Fq}^2$ if $\parameter = 0$.

In this article, we prove this conjecture by computing explicitly $\card{\curveA{\parameter}(\Fqr)}$ and $\card{\curveB{\parameter}(\Fqr)}$ in terms of Gauss sums using formulas given by Koblitz in \cite[\S 5]{KOB.nb.points.hyper} and comparing them to those given for $\card{\quintic{\parameter}(\Fqr)}$ by Candelas, de la Ossa, and Rodriguez Villegas. More precisely, we show the following theorem (with analogous formulas for $\curveB{\parameter}$).

\begin{thm}\label{main.theorem}
If $\curveA{\parameter}$ is the affine curve $y^5 = x^2(1-x)^3(x-\parameter^5)^2$,
\begin{align}\label{aim:Na=NA:psi<>0}&\card{\curveA{\parameter}(\Fqr)} = q^r - N_a(q^r) &\quad \text{if $\parameter \neq 0$} \\
\label{aim:Na=NA:psi=0}&\card{\curveA{\parameter}(\Fqr)} = q^r - \frac{1}{2}N_a(q^r) &\quad \text{if $\parameter = 0$}\end{align}
\end{thm}

There are a number of remarks that can be made about the two factors $P_{\curveA{\parameter}}$ and $P_{\curveB{\parameter}}$.

\begin{rem}\label{remark:PA.PB:rho-th.power} If $\rho \in \{1,2,4\}$ is the order of $q$ in $(\Zset/5\Zset)^\times$, then $P_{\curveA{\parameter}/\Fq}(t) = P_{\curveA{\parameter}/\Fqrho}(t^\rho)^{1/\rho}$ which implies that $P_{\curveA{\parameter}/\Fq}(t)$ is the $\rho$-th power of a polynomial. The same is true for $\curveB{\parameter}$.
\end{rem}

\begin{rem}\label{remark:PA.PB.psi<>0:squares} If $\parameter\neq 0$ and $\parameter^5 \neq 1$, then $P_{\curveA{\parameter}/\Fq}$ and $P_{\curveB{\parameter}/\Fq}$ are squares.
\end{rem}

\begin{rem}\label{remark:PA.PB.psi<>0:4-th.power} If $\parameter = 0$, $p \not\equiv 1 \mod 5$ and $q \equiv 1 \mod 5$, then $P_{\curveA{\parameter}/\Fq}(t) = P_{\curveB{\parameter}/\Fq}(t) = (1-q't)^4$ where $q'$ is the square root of $q$ which is $\equiv 1 \mod 5$.\end{rem}

By combining remark~\ref{remark:PA.PB:rho-th.power} and remark~\ref{remark:PA.PB.psi<>0:4-th.power}, one can easily prove the observations of \cite[\S 12 p.~132]{CdlORV.II} when $\parameter = 0$ and $\rho = 2$ or $\rho = 4$. For a proof of remark~\ref{remark:PA.PB.psi<>0:squares}, we refer to \cite[\S 11.1 p.~129-130]{CdlORV.II}; the argument is that it is possible to transform $\curveA{\parameter}$ and $\curveB{\parameter}$ into hyperelliptic curves and then use the existence of a special automorphism of these hyperelliptic curves to deduce that the Jacobian of these two curves is isogenous to a square.

The article is organized as follows. In section~\ref{section:gauss.sums}, we recall all the formulas on Gauss and Jacobi sums we will need. In section~\ref{section:nb.points.curves} we give the formulas for $\card{\curveA{\parameter}(\Fqr)}$, $\card{\curveB{\parameter}(\Fqr)}$ and prove remarks~\ref{remark:PA.PB:rho-th.power} and \ref{remark:PA.PB.psi<>0:4-th.power}. In section~\ref{section:diagonal.case} and section~\ref{section:generic.case}, we prove Theorem~\ref{main.theorem} when $\parameter = 0$ and $\parameter \neq 0$ respectively, and finally, in section~\ref{section:singular.case}, we mention what happens when the quintic is singular.

The method we use does not give any information on a geometric link between the two curves $\curveA{\parameter}$, $\curveB{\parameter}$ and the quintic $\quintic{\parameter}$, but can be generalized to give the explicit factorization of the zeta function of $x_1^n + \dots + x_n^n - n \parameter x_1 \dots x_n = 0$, at least when $n$ is a prime number. This will be the subject of a subsequent article.

\section{Gauss and Jacobi sums formulas}\label{section:gauss.sums}

\begin{defn}[Gauss sums]
Let $\bbOmega$ be an algebraically closed field of characteristic zero. Let $\caradd : \Fq \to \bbOmega^*$ be a non trivial additive character of $\Fq$. For any character $\carmult : \Fq^* \to \bbOmega^*$ of $\Fq^*$, define the Gauss sum $\gausssum{\carmult}$ as being
\[\gausssum{\carmult} = \sum_{x\in\Fq^*}{\caradd(x)\carmult(x)}.\]
\end{defn}

If $\carmulttriv$ is the trivial character of $\Fq^*$, we have $\gausssum{\carmulttriv} = -1$. Note that $\gausssum{\carmult^i}$ only depends on $i \mod q-1$.

\begin{prop}[Reflection formula]
If $\carmult$ is a non trivial character of $\Fq^*$,
\begin{equation}\label{formula:relection}
\gausssum{\carmult} \gausssum{\carmult^{-1}} = \carmult(-1) q.
\end{equation}
\end{prop}

\begin{pf}
See \cite[Theorem~1.1.4~(a) p.~10]{BEW}. \qed
\end{pf}

\begin{prop}[Hasse-Davenport formula]
If $\carmult$ is a character of $\Fq^*$,
\begin{equation}\label{formula:hasse-davenport}
-\gausssumadd{\caradd \circ \TrFqrFq}{\carmult \circ \NFqrFq} = (-\gausssum{\carmult})^r.
\end{equation}
\end{prop}

\begin{pf}
See \cite[Theorem~11.5.2 p.~360]{BEW}. \qed
\end{pf}

\begin{prop}[Purity formula]
Let $d \geq 3$ be an integer and $\caraddbis$ a non trivial additive character of $\Fp$. Suppose there exists an integer $s$ such that $p^s \equiv -1 \mod d$. If $\sigma$ is the smallest integer such that $p^\sigma \equiv -1 \mod d$ and if we let $q = p^{2 \sigma m}$ (so that $\sqrt{q} = p^{\sigma m}$ is an integer), then, for any multiplicative character $\carmult$ of order $d$ of $\Fq^*$,
\begin{equation}\label{formula:purity}
\gausssumadd{\caraddbis \circ \TrFqFp}{\carmult} = \begin{cases} (-1)^{m-1} \sqrt{q} & \text{if $d$ is odd or if $\frac{p^\sigma+1}{d}$ is even} \\ -\sqrt{q} & \text{if $d$ is even and if $\frac{p^\sigma+1}{d}$ is odd}\end{cases}
\end{equation}
\end{prop}

\begin{pf}
See \cite[Theorem~11.6.3 p.~364]{BEW}. \qed
\end{pf}

\begin{prop}[Multiplication formula]
Let $d \geq 1$ be a divisor of $q-1$. If $\carmult$ is a character of $\Fq^*$,
\begin{equation}\label{formula:multiplication}
\frac{\gausssum{\carmult^d}}{\prod_{\carmult '^d = \carmulttriv}{\gausssum{\carmult\carmult '}}} = \frac{\carmult(d)^d}{\prod_{\substack{\carmult '^d = \carmulttriv \\ \carmult ' \neq \carmulttriv}}{\gausssum{\carmult '}}}.
\end{equation}
\end{prop}

\begin{pf}
See \cite[Theorem~11.3.5 p.~355]{BEW}. \qed
\end{pf}

\begin{rem}
The product $\prod_{\substack{\carmult '^d = \carmulttriv \\ \carmult ' \neq \carmulttriv}}{\gausssum{\carmult '}}$ is given by
\begin{equation}\label{formula:product.simplification}
\begin{cases} q^{\frac{d-1}{2}} & \text{if $d$ is odd} \\ (-1)^{\frac{(q-1)(d-2)}{8}} q^{\frac{d-2}{2}} \gausssum{\carmultquad} & \text{if $d$ is even}\end{cases}
\end{equation}
where $\carmultquad$ is the character of order $2$ of $\Fq^*$. This formula is a direct consequence of the reflection formula~(\ref{formula:relection}) by grouping Gauss sums by pair. When $\bbOmega = \Cset$, it is possible to give an exact formula for $\gausssum{\carmultquad}$ (see \cite[Theorem~11.5.4 p.~362]{BEW}), but we will not have any use of it as we will only be concerned with the case $d=5$.
\end{rem}

\begin{defn}[Jacobi sums]
If $\carmult$ and $\carmult '$ are two characters of $\Fq^*$, define the corresponding Jacobi sum as being
\[\jacobisum(\carmult,\carmult ') = \sum_{\substack{x,x '\in\Fq^*\\x+x '=1}} {\carmult(x) \carmult '(x ')} = \sum_{\substack{x\in\Fq\\x\neq 0, x\neq 1}} {\carmult(x) \carmult '(1-x)}.\]
\end{defn}

\begin{prop}[Link with Gauss sums]
If $\carmult$ and $\carmult '$ are two characters of $\Fq^*$,
\begin{equation}\label{formula:jacobi-gauss}
\jacobisum(\carmult,\carmult ') = \begin{cases} ~\frac{\gausssum{\carmult} \gausssum{\carmult '}}{ \gausssum{\carmult \carmult '}} & \text{if $\carmult\carmult ' \neq \carmulttriv$} \\ \frac{1}{q} \frac{\gausssum{\carmult} \gausssum{\carmult '}}{ \gausssum{\carmult \carmult '}} & \text{if $\carmult\carmult ' = \carmulttriv$ and $\carmult \neq \carmulttriv$}\end{cases}
\end{equation}
\end{prop}

\begin{pf}
For the first formula, see \cite[Eq.~(11.6.4) p.~365]{BEW} (notice that there is nothing to prove if $\carmult$ or $\carmult '$ is trivial). For the second formula, use \cite[Eq.~(11.6.3) p.~365]{BEW} and then the reflection formula~(\ref{formula:relection}) to see that $J(\carmult,\carmult^{-1}) = -\carmult(-1) = \frac{1}{q} \frac{G(\caradd,\carmult) G(\caradd,\carmult^{-1})}{ G(\caradd,\carmult\carmult^{-1})}.$ \qed
\end{pf}

\begin{prop}[Jacobi sum of inverses]
If $\carmult$ is a non trivial character of $\Fq^*$ such that $\carmult(-1)=1$ and if $\carmultbis$ is an arbitrary character of $\Fq^*$,
\begin{equation}\label{formula:jacobisum(inv,inv)}\jacobisum(\carmult^{-1},\carmultbis^{-1}) = \begin{cases} \frac{1}{q} \frac{\gausssum{\carmult^{-1}}\gausssum{\carmult\carmultbis}}{\gausssum{\carmultbis}} & \text{if $\carmultbis = \carmulttriv$} \\ ~~\frac{\gausssum{\carmult^{-1}}\gausssum{\carmult\carmultbis}}{\gausssum{\carmultbis}} & \text{otherwise}\end{cases}
\end{equation}
\end{prop}

\begin{pf}
If $\carmult$ and $\carmultbis$ are two characters of $\Fq^*$, we have
\[\begin{split}
\jacobisum(\carmult^{-1},\carmult\carmultbis)
& = \sum_{\substack{x,y\in\Fq^* \\ x+y=1}}{\carmult^{-1}(x)\carmult\carmultbis(y)} = \sum_{\substack{x,y\in\Fq^* \\ x+y=1}}{\carmult^{-1}(\tfrac{x}{y})\carmultbis^{-1}(\tfrac{1}{y})}~; \\
& = \sum_{\substack{x',y'\in\Fq^* \\ x'+y'=1}}{\carmult^{-1}(-x')\carmultbis^{-1}(y')}  = \carmult(-1)\jacobisum(\carmult^{-1},\carmultbis^{-1}),
\end{split}\]
where we made the change of variables $x' = -\tfrac{x}{y}$ et $y' = \tfrac{1}{y}$. Combining this formula with formula~(\ref{formula:jacobi-gauss}), we get at once formula~(\ref{formula:jacobisum(inv,inv)}) because $\carmult$ is non trivial and $\carmult(-1)=1$. \qed
\end{pf}

\begin{prop}[Fourier inversion formula]
If $f : \Fq^* \to \bbOmega$ is a map,
\begin{equation}\label{formula:Fourier.inversion}
\forall \parametercurve \in \Fq^*, \quad f(\parametercurve) = \frac{1}{q-1}\sum_{\carmultbis \in \Fqhat}{\left(\sum_{\parametercurvebis \in \Fq^*}{f(\parametercurvebis)\carmultbis^{-1}(\parametercurvebis)}\right) \carmultbis(\parametercurve)}.
\end{equation}
\end{prop}

\begin{pf}
It is a direct consequence of the orthogonality relations for characters of the finite abelian group $\Fq^*$. \qed
\end{pf}

\begin{rem}\label{remark:choice.caradd}
From now on, we will take for additive character $\caradd$ a character of the form $\caraddbis \circ \TrFqFp$ where $\caraddbis$ is a non trivial additive character of $\Fp$, so purity formula~(\ref{formula:purity}) will be valid.
\end{rem}

\section{Number of points of the hypergeometric curves}\label{section:nb.points.curves}

In the introduction, we gave equations $y^5 = x^2(1-x)^3(x-\parameter^5)^2$ and $y^5 = x^2(1-x)^4(x-\parameter^5)$ for $\curveA{\parameter}$ and $\curveB{\parameter}$ respectively. When $\parameter \neq 0$, we transform these equations into
\[y^5 = x^2(1-x)^3(1-\tfrac{1}{\parameter^5}x)^2 \quad \text{ and } \quad y^5 = x^2(1-x)^4(1-\tfrac{1}{\parameter^5}x),\]
which does not change the number of points and is more convenient for our computations.

\subsection{General remarks}\label{section:nb.points.curves:general.remarks}

To show Theorem~\ref{main.theorem}, we only have to consider the case $r=1$ as $q$ is arbitrary. Moreover, the following lemma shows that we only need to consider the case $q \equiv 1 \mod 5$ to compute $\card{\curveA{\parameter}(\Fq)}$.

\begin{lem}\label{result:nb.zerosA.trivial}
Let $Q$ be a polynomial with coefficients in $\Fq$ and $\curveC$ the affine curve $y^5 = Q(x)$. If $q \not\equiv 1 \mod 5$, then $\card{\curveC(\Fq)} = q$.
\end{lem}

\begin{pf}
Because $q \not\equiv 1 \mod 5$, the map $y \mapsto y^5$ is a bijection of $\Fq$ onto itself and so $\curveC$ has the same number of points than the curve $y = Q(x)$ i.e. $q$ points. \qed
\end{pf}

\begin{rem}\label{remark:PA.PB:rho-th.power:proof}
This lemma shows remark~\ref{remark:PA.PB:rho-th.power} of the introduction.
\end{rem}

\subsection{The case $\parameter = 0$}\label{section:nb.points.curves:parameter=0}

The curves $\curveA{0}$ and $\curveB{0}$ are isomorphic so have the same number of points. It is thus sufficient to give a formula for $\card{\curveA{0}(\Fq)}$.

\begin{thm}\label{result:nbzeros.AB.psi=0}
If $q \equiv 1 \mod 5$ and if $\curveA{0}$ is the affine curve $y^5 = x^4(1-x)^3$,
\[\card{\curveA{0}(\Fq)} = q - \sum_{\substack{\carmult^5 = \carmulttriv \\ \carmult \neq \carmulttriv}}{\left(-\frac{1}{q}\gausssum{\carmult}^2 \gausssum{\carmult^3}\right)}.\]
\end{thm}

\begin{pf}
We recall the proof of this classical result (see for example \cite[\S 1, p.~202-203]{GR} or \cite[\S 5, p.~20]{KOB.nb.points.hyper}). It is a straightforward application of the method Weil used in \cite{WEIL.nb.sol.finite.fields}.

The equation of $\curveA{0}$ is of the type $y^5 = Q(x)$ where $Q$ is a polynomial. The number of points of the affine curve $\curveA{0}$ is given by
\[\card{\curveA{0}(\Fq)} = \card{\set{(x,y) \in \Fq \times \Fq}{y^5 = Q(x)}}.\]
The fundamental remark is that, if $z \in \Fq$,
\[\card{\set{y \in \Fq}{y^5 = z}} = \begin{cases} 1 & \text{if $z=0$} \\ 1 + \sum_{\substack{\carmult^5 = \carmulttriv \\ \carmult \neq \carmulttriv}}{\carmult(z)} & \text{if $z \neq 0$}\end{cases}\]
Thus,
\[\begin{split}
\card{\curveA{0}(\Fq)}
& = \card{\set{(x,y) \in \Fq \times \Fq}{y^5 = Q(x) \text{ and } Q(x) = 0}} \\
& \quad + \card{\set{(x,y) \in \Fq \times \Fq}{y^5 = Q(x) \text{ and } Q(x) \neq 0}} ;\\
& = \sum_{\substack{x \in \Fq \\ Q(x) = 0}}{1} + \sum_{\substack{x \in \Fq \\ Q(x) \neq 0}}{\Bigg(1+\sum_{\substack{\carmult^5 = \carmulttriv \\ \carmult \neq \carmulttriv}}{\carmult(Q(x))}\Bigg)},
\end{split}\]
and so
\begin{equation}\label{formula:nb.zeros.curve}
\card{\curveA{0}(\Fq)} = q + \sum_{\substack{\carmult^5 = \carmulttriv \\ \carmult \neq \carmulttriv}}{\sum_{\substack{x \in \Fq \\ Q(x) \neq 0}}{\carmult(Q(x))}}.
\end{equation}
Replacing $Q(x)$ by $x^4(1-x)^3$, we obtain
\[\begin{split}
\card{\curveA{0}(\Fq)}
& = q + \sum_{\substack{\carmult^5 = \carmulttriv \\ \carmult \neq \carmulttriv}}{\sum_{\substack{x \in \Fq \\ x \neq 0 , x \neq 1}}{\carmult^4(x) \carmult^3(1-x)}} = q + \sum_{\substack{\carmult^5 = \carmulttriv \\ \carmult \neq \carmulttriv}}{\jacobisum(\carmult^4,\carmult^3)}.
\end{split}\]
We now use formula~(\ref{formula:jacobi-gauss}) to express the Jacobi sums in terms of Gauss sums
\[\begin{split}
\card{\curveA{0}(\Fq)}
& = q + \sum_{\substack{\carmult^5 = \carmulttriv \\ \carmult \neq \carmulttriv}}{\frac{\gausssum{\carmult^4}\gausssum{\carmult^3}} {\gausssum{\carmult^7}}} = q + \sum_{\substack{\carmult^5 = \carmulttriv \\ \carmult \neq \carmulttriv}}{\frac{\gausssum{\carmult^4}\gausssum{\carmult^3}} {\gausssum{\carmult^2}}} ; \\
& = q + \sum_{\substack{\carmult^5 = \carmulttriv \\ \carmult \neq \carmulttriv}}{\frac{1}{q} \gausssum{\carmult^4}\gausssum{\carmult^3}^2 },
\end{split}\]
by using the reflection formula~(\ref{formula:relection}) and the fact that $\carmult(-1) = 1$ (because $\carmult^5 = \carmulttriv$). To obtain the announced formula, we just make the change of variable $\carmult \mapsto \carmult^2$. \qed
\end{pf}

\begin{rem}
Assume as before that $q \equiv 1 \mod 5$ and replace $\Fq$ by $\Fqr$ in the formula for $\card{\curveA{\parameter}(\Fq)}$ we have just obtained. As the characters of order $5$ of $\Fqr^*$ are exactly the $\carmult \circ \TrFqrFq$ where $\carmult$ is a character of order $5$ of $\Fq^*$, the choice of additive character made in remark~\ref{remark:choice.caradd} and the Hasse-Davenport formula~(\ref{formula:hasse-davenport}) show that
\[\card{\curveA{0}(\Fqr)} = q^r - \sum_{\substack{\carmult \in \Fqhat \\ \carmult^5 = \carmulttriv, \carmult \neq \carmulttriv}}{\left(-\frac{1}{q} \gausssum{\carmult}^2\gausssum{\carmult^3}\right)^r},\]
and so, when $q \equiv 1 \mod 5$, the zeta function is given by
\[\zetafunction{\curveA{0}/\Fq}{t} = \frac{\prod_{\substack{\carmult^5 = \carmulttriv \\ \carmult \neq \carmulttriv}}{\left(1+\frac{1}{q} \gausssum{\carmult}^2\gausssum{\carmult^3} t\right)}}{1-qt}.\]
\end{rem}

\begin{rem}\label{remark:PA.PB.psi<>0:4-th.power:proof}
Let us deduce remark~\ref{remark:PA.PB.psi<>0:4-th.power} of the introduction. Because $p \not\equiv 1 \mod 5$ and $q \equiv 1 \mod 5$, we must have $q = p^{2m}$ if $p \equiv -1 \mod 5$ and $q = p^{4m}$ if $p \equiv \pm 2 \mod 5$. The purity formula~(\ref{formula:purity}) then shows that the numerator of the zeta function of $\curveA{0}$ is
\[(1-(-1)^m \sqrt{q} t)^4 = (1-q't)^4,\]
where $q' = (-1)^m \sqrt{q}$ is the square root of $q$ which is $\equiv 1 \mod 5$.
\end{rem}

\subsection{The case $\parameter \neq 0$}

When $\parameter \neq 0$, the equations of both $\curveA{\parameter}$ and $\curveB{\parameter}$ are of the type $\curveC[\parametercurve]$~:~$y^5 = x^\expcurvea (1-x)^\expcurveb (1-\parametercurve x)^{5-\expcurveb}$ where $\parametercurve = \tfrac{1}{\parameter^5}$ is a fifth power.

\begin{thm}\label{result:nbzeros.AB.psi<>0}
If $q \equiv 1 \mod 5$ and if $\curveC[\parametercurve]$ is the affine curve $y^5 = x^\expcurvea (1-x)^\expcurveb (1-\parametercurve x)^{5-\expcurveb}$ with $\parametercurve \in \Fq^*$ and $1 \leq \expcurvea , \expcurveb\leq 4$,
\[\card{\curveC[\parametercurve](\Fq)} = q + \sum_{\substack{\carmult^5 = \carmulttriv \\ \carmult \neq \carmulttriv}}{\left(\frac{1}{q-1} \sum_{\carmultbis \in \Fqhat}{q^{1-\nu} \frac{\gausssum{\carmult^\expcurvea\carmultbis} \gausssum{\carmult^\expcurveb\carmultbis}}{\gausssum{\carmultbis} \gausssum{\carmult^{\expcurvea+\expcurveb}\carmultbis}} \carmultbis(\parametercurve)}\right)},\]
where $\nu$ is the number of trivial characters in the pair $(\carmultbis,\carmult^{\expcurvea+\expcurveb}\carmultbis)$.
\end{thm}

\begin{pf}
We follow the proof of Koblitz \cite[Theorem~3 p.~18]{KOB.nb.points.hyper} and then derive the above formula.

Formula~(\ref{formula:nb.zeros.curve}) is valid for $\curveC[\parametercurve]$ and gives
\[\begin{split}
\card{\curveC[\parametercurve](\Fq)}
& = q + \sum_{\substack{\carmult^5 = \carmulttriv\\\carmult \neq \carmulttriv}}{\sum_{\substack{x \in \Fq \\ x \neq 0, x \neq 1, x \neq 1/\parametercurve}} {\carmult^{\expcurvea}(x) \carmult^{\expcurveb}(1-x) \carmult^{-\expcurveb}(1-\parametercurve x)}} ;\\
& = q + \sum_{\substack{\carmult^5 = \carmulttriv\\\carmult \neq \carmulttriv}}{\nbzeros{\curveC[\parametercurve]/\Fq,\carmult}},
\end{split}\]
where, if $\carmult$ is a non trivial character of order $5$ of $\Fq^*$,
\[\nbzeros{\curveC[\parametercurve]/\Fq,\carmult} = \sum_{\substack{x \in \Fq \\ x \neq 0, x \neq 1, x \neq 1/\parametercurve}} {\carmult^{\expcurvea}(x) \carmult^{\expcurveb}(1-x) \carmult^{-\expcurveb}(1-\parametercurve x)}.\]
Following Koblitz \cite[p.~19]{KOB.nb.points.hyper}, we make a Fourier transform. We choose a character $\carmultbis$ of $\Fq^*$, multiply $\nbzeros{\curveC[\parametercurve]/\Fq,\carmult}$ by $\carmultbis^{-1}(\parametercurve)$, and sum on all $\parametercurve \neq 0$
\[\sum_{\parametercurve\in\Fq^*}{\nbzeros{\curveC[\parametercurve]/\Fq,\carmult}\carmultbis^{-1}(\parametercurve)} = \sum_{\substack{x \in \Fq,\parametercurve\in\Fq^* \\ x \neq 0, x \neq 1, x \neq 1/\parametercurve}} {\carmult^{\expcurvea}(x) \carmult^{\expcurveb}(1-x) \carmult^{-\expcurveb}(1-\parametercurve x) \carmultbis^{-1}(\parametercurve)}.\]
We now make the change of variable $t = \parametercurve x$
\[\sum_{\parametercurve\in\Fq^*}{\nbzeros{\curveC[\parametercurve]/\Fq,\carmult}\carmultbis^{-1}(\parametercurve)} = \sum_{\substack{x \in \Fq,t\in\Fq\\x\neq 0,x\neq 1,t\neq 0, t \neq 1}} {\carmult^{\expcurvea}(x) \carmult^{\expcurveb}(1-x) \carmult^{-\expcurveb}(1-t) \carmultbis^{-1}(t) \carmultbis(x)}.\]
Grouping the terms with only $x$ together and only $t$ together, we obtain two Jacobi sums
\[\begin{split}
\sum_{\parametercurve\in\Fq^*}{\nbzeros{\curveC[\parametercurve]/\Fq,\carmult}\carmultbis^{-1}(\parametercurve)} & = \Bigg(\sum_{\substack{x \in \Fq\\x\neq 0, x \neq 1}} {(\carmult^{\expcurvea}\carmultbis)(x) \carmult^{\expcurveb}(1-x)}\Bigg) \Bigg(\sum_{\substack{t \in \Fq \\ t\neq 0, t \neq 1 }}{\carmult^{-\expcurveb}(1-t) \carmultbis^{-1}(t)}\Bigg);\\
& = \jacobisum(\carmult^{\expcurvea}\carmultbis,\carmult^{\expcurveb}) \jacobisum(\carmult^{-\expcurveb},\carmultbis^{-1}).
\end{split}\]
We now use the Fourier inversion formula~(\ref{formula:Fourier.inversion})
\[\begin{split}
\nbzeros{\curveC[\parametercurve]/\Fq,\carmult}
& = \frac{1}{q-1}\sum_{\carmultbis \in \Fqhat}{\left(\sum_{\parametercurvebis\in\Fq^*} {\nbzeros{\curveC[\parametercurvebis]/\Fq,\carmult} \carmultbis^{-1}(\parametercurvebis)}\right) \carmultbis(\parametercurve)};\\
& = \frac{1}{q-1}\sum_{\carmultbis \in \Fqhat}{\jacobisum(\carmult^{\expcurvea}\carmultbis,\carmult^{\expcurveb}) \jacobisum(\carmult^{-\expcurveb},\carmultbis^{-1}) \carmultbis(\parametercurve)}.
\end{split}\]
Because $\carmult^5 = \carmulttriv$, we have $\carmult(-1) = 1$ and so we can use formulas~(\ref{formula:jacobi-gauss}) and (\ref{formula:jacobisum(inv,inv)}) to express these Jacobi sums in terms of Gauss sums and obtain
\[\jacobisum(\carmult^{\expcurvea}\carmultbis,\carmult^{\expcurveb}) \jacobisum(\carmult^{-\expcurveb},\carmultbis^{-1}) = \frac{1}{q^{\nu}} \frac{\gausssum{\carmult^\expcurvea\carmultbis} \gausssum{\carmult^\expcurveb}} {\gausssum{\carmult^{\expcurvea+\expcurveb}\carmultbis}} \frac{\gausssum{\carmult^\expcurveb\carmultbis} \gausssum{\carmult^{-\expcurveb}}}{\gausssum{\carmultbis}}.\]
Finally, using the reflection formula~(\ref{formula:relection}),
\begin{equation}\label{result:nbzeros.AB.carmult}
\nbzeros{\curveC[\parametercurve]/\Fq,\carmult} = \frac{1}{q-1} \sum_{\carmultbis \in \Fqhat}{q^{1-\nu} \frac{\gausssum{\carmult^\expcurvea\carmultbis} \gausssum{\carmult^\expcurveb\carmultbis}}{\gausssum{\carmultbis} \gausssum{\carmult^{\expcurvea+\expcurveb}\carmultbis}} \carmultbis(\parametercurve)},
\end{equation}
which gives the announced formula. \qed
\end{pf}

We now rewrite the formula for $\card{\curveC[\parametercurve](\Fq)}$ of the previous theorem in a form better suited for our use of it in \S\ref{section:generic.case}.

\begin{cor}\label{result:nbzeros.AB:decomp}
Suppose that $q \equiv 1 \mod 5$ and let $\curveC[\parametercurve]$ be the affine curve $y^5 = {x^\expcurvea (1-x)^\expcurveb (1-\parametercurve x)^{5-\expcurveb}}$ with $\parametercurve \in \Fq^*$ a fifth power and $1 \leq \expcurvea , \expcurveb\leq 4$. If $\carmult$ is a non trivial character of order $5$ of $\Fq^*$, we can write
\[\card{\curveC[\parametercurve](\Fq)} = q + 2 \nbzeros{\curveC[\parametercurve]/\Fq,\carmult} + 2 \nbzeros{\curveC[\parametercurve]/\Fq,\carmult^2},\]
where $\nbzeros{\curveC[\parametercurve]/\Fq,\carmult}$ is given by~(\ref{result:nbzeros.AB.carmult}).
\end{cor}

\begin{pf}
It is a simple consequence of Theorem~\ref{result:nbzeros.AB.psi<>0} and of the formula $\nbzeros{\curveC[\parametercurve]/\Fq,\carmult '} = \nbzeros{\curveC[\parametercurve]/\Fq,\carmult '^{-1}}$. To prove this last formula, write
\[\nbzeros{\curveC[\parametercurve]/\Fq,\carmult '^{-1}} = \frac{1}{q-1} \sum_{\carmultbis \in \Fqhat}{q^{1-\nu} \frac{\gausssum{\carmult '^{-\expcurvea}\carmultbis} \gausssum{\carmult '^{-\expcurveb}\carmultbis}}{\gausssum{\carmultbis} \gausssum{\carmult '^{-(\expcurvea+\expcurveb)}\carmultbis}} \carmultbis(\parametercurve)},\]
and make the change of variable $\carmultbis = \carmult '^{a+b} \carmultbis '$. It transforms $\nbzeros{\curveC[\parametercurve]/\Fq,\carmult '^{-1}}$ into $\nbzeros{\curveC[\parametercurve]/\Fq,\carmult '}$ because $\parametercurve$ is a fifth power and because the number of trivial characters in the pair $(\carmultbis,\carmult '^{-(\expcurvea+\expcurveb)}\carmultbis)$ in the same than in the pair $(\carmultbis ',{\carmult '}^{\expcurvea+\expcurveb}\carmultbis ')$. \qed
\end{pf}

\section{Number of points of the quintic in the diagonal case}\label{section:diagonal.case}

The aim of this section is to prove formula~(\ref{aim:Na=NA:psi=0}) of Theorem~\ref{main.theorem} of the introduction. As in the previous section, we restrict ourselves to the case $r = 1$ as $q$ is arbitrary. We begin with the case $q \not\equiv 1 \mod 5$.

\begin{lem}\label{result:nb.zeros:quintic:trivial}
If $q \not\equiv 1 \mod 5$, the number of points of $\quintic{0} : x_1^5 + \dots + x_5^5 = 0$ in $\Proj^4_{\Fq}$ is given by
\[\card{\quintic{0}(\Fq)} = 1 + q + q^2 + q^3.\]
\end{lem}

\begin{pf}
When $q \not\equiv 1 \mod 5$, the map $x \mapsto x^5$ is a bijection of $\Fq$ onto itself, so $\quintic{0}$ has the same number of points as the hyperplane $x_1 + \dots + x_5 = 0$. \qed
\end{pf}

Comparing with Lemma~\ref{result:nb.zerosA.trivial}, this proves formula~(\ref{aim:Na=NA:psi=0}) of Theorem~\ref{main.theorem} of the introduction when $q \not\equiv 1 \mod 5$. We now proceed to the case $q \equiv 1 \mod 5$.

\begin{thm}
If $q \equiv 1 \mod 5$, we can write
\[\card{\quintic{0}(\Fq)} = 1 + q + q^2 + q^3 - N_w(q) - 25 q N_a(q),\]
where
\[\begin{split}
& N_w(q) = \sum_{\substack{\carmult^5 = \carmulttriv \\ \carmult \neq \carmulttriv}}{\left(-\frac{1}{q} \gausssum{\carmult}^5\right)};\\
& N_a(q) = 2\sum_{\substack{\carmult^5 = \carmulttriv \\ \carmult \neq \carmulttriv}}{\left(-\frac{1}{q}\gausssum{\carmult}^2\gausssum{\carmult^3}\right)}.
\end{split}\]
\end{thm}

\begin{pf}
This is a direct consequence of Weil's formula \cite[Eq.~(3) p.~500]{WEIL.nb.sol.finite.fields}. More precisely, if we chose a non trivial character $\carmult$ of order $5$ of $\Fq^*$, we have
\[\card{\quintic{0}(\Fq)} = 1 + q + q^2 + q^3 - \sum_{\substack{1 \leq s_1,\dots,s_5 \leq 4 \\ s_1 + \dots + s_5 \equiv 0 \mod 5}}{\left(-\frac{1}{q} \gausssum{\carmult^{s_1}} \dots \gausssum{\carmult^{s_5}}\right)}.\]
Up to permutation, there are only twelve $5$-uples $(s_1,\dots,s_5)$ that index this sum. Explicitly
\begin{center}\begin{tabular}{cc}
$(s_1,\dots,s_5)$ & \textsc{\# permutations} \\
$(1,1,1,1,1),(2,2,2,2,2),(3,3,3,3,3),(4,4,4,4,4)$ & $ 1$ \\
$(1,1,3,1,4),(2,2,1,2,3),(3,3,4,2,3),(4,4,2,1,4)$ & $20$ \\
$(1,1,3,2,3),(2,2,1,1,4),(3,3,4,1,4),(4,4,2,2,3)$ & $30$ \\
\end{tabular}\end{center}
We are thus able to enumerate all the products of Gauss sums appearing in the previous formula for $\card{\quintic{0}(\Fq)}$. Using the reflection formula~(\ref{formula:relection}) and grouping together the terms, we get the following table.
\begin{center}\begin{tabular}{ccc}
$(s_1,\dots,s_5)$ & $\frac{1}{q} G(\caradd,\carmult^{s_1}) \dots G(\caradd,\carmult^{s_5})$ & \textsc{multiplicity}\\
$(j,j,j,j,j)$ & $\frac{1}{q} G(\caradd,\carmult^j)^5$ & $1$\\
$(1,1,3,1,4)$ and $(1,1,3,2,3)$ & $G(\caradd,\carmult)^2 G(\caradd,\carmult^3)$ & $20+30=50$\\
$(2,2,1,2,3)$ and $(2,2,1,1,4)$ & $G(\caradd,\carmult^2)^2 G(\caradd,\carmult)$ & $20+30=50$\\
$(3,3,4,2,3)$ and $(3,3,4,1,4)$ & $G(\caradd,\carmult^3)^2 G(\caradd,\carmult^4)$ & $20+30=50$\\
$(4,4,2,1,4)$ and $(4,4,2,2,3)$ & $G(\caradd,\carmult^4)^2 G(\caradd,\carmult^2)$ & $20+30=50$\\
\end{tabular}\end{center}
This gives the formula we want, namely
\begin{multline*}
\card{\quintic{0}(\Fq)} = 1 + q + q^2 + q^3 - \sum_{j=1}^{4}{\left(-\frac{1}{q} \gausssum{\carmult^j}^5\right)} \\ - 50 \sum_{j=1}^{4}{\left(- \gausssum{\carmult^j}^2 \gausssum{\carmult^{3j}}\right)}. \qed
\end{multline*}
\end{pf}

By comparing with Theorem~\ref{result:nbzeros.AB.psi=0}, we immediately obtain formula~(\ref{aim:Na=NA:psi=0}) of Theorem~\ref{main.theorem} of the introduction when $q \equiv 1 \mod 5$.

\begin{rem}\label{remark:H0}
When $\parameter = 0$, the factor $N_w(q)$ comes from the hypergeometric hypersurface
\[H_0~:~y^5 = x_1 x_2 x_3 (1-x_1-x_2-x_3).\]
More precisely, using the techniques from sections~\ref{section:nb.points.curves:general.remarks} and \ref{section:nb.points.curves:parameter=0}, it is straightforward to show that
\[\card{H_0(\Fq)} = q^3 - N_w(q).\]
This is in accordance with \cite[Proposition~3.2]{GY}. 
\end{rem}

\section{Number of points of the quintic in the generic case}\label{section:generic.case}

We now consider the case $\parameter \neq 0$ and proceed to show formula~(\ref{aim:Na=NA:psi<>0}) of Theorem~\ref{main.theorem} of the introduction. Let us first recall the result shown by Candelas, de la Ossa, and Rodriguez Villegas. Just like in the two previous sections, we restrict ourselves to the case $r=1$ as $q$ is arbitrary.

\begin{thm}[Candelas, de la Ossa, Rodriguez Villegas]\label{result:CdlORV.psi<>0}
If $\parameter \neq 0$, then
\[\card{\quintic{\parameter}(\Fq)} = 1 + q + q^2 + q^3 - N_w(q) - 10 q N_a(q) - 15 q N_b(q) + 24 N_\textup{sing}(q),\]
where $N_a(q)$, $N_b(q)$ and $N_\textup{sing}(q)$ vanish if $q \not\equiv 1 \mod 5$. Moreover,
\[\card{\mirrorquintic{\parameter}(\Fq)} = 1 + q + q^2 + q^3 - N_w(q) + 100(q+q^2),\]
where $\mirrorquintic{\parameter}$ is the ``mirror'' of $\quintic{\parameter}$. We also have the following formulas for $N_a(q)$, $N_b(q)$ and $N_\textup{sing}(q)$ when $q \equiv 1 \mod 5$
\[\begin{split}
 -N_{a}(q) &= \tfrac{2}{q}N_{(0,0,0,1,4),\carmult}(q) + \tfrac{2}{q}N_{(0,0,0,2,3),\carmult}(q);\\
 -N_{b}(q) &= \tfrac{2}{q}N_{(0,0,1,1,3),\carmult}(q) + \tfrac{2}{q}N_{(0,0,1,2,2),\carmult}(q);\\
 N_{\textup{sing}}(q) &= N_{(0,1,2,3,4),\carmult}(q),
\end{split}\]
where $\carmult$ is a fixed character of order $5$ of $\Fq^*$ and
\[N_{(s_1,s_2,s_3,s_4,s_5),\carmult}(q) = \frac{1}{q-1} \sum_{\carmultbis \in \Fqhat}{q^{5-z-\delta}\frac{\gausssum{\carmultbis^{5}}} {\prod_{j=1}^{5}{\gausssum{\carmult^{s_j}\carmultbis}}} \carmultbis(\tfrac{1}{(5\parameter)^5})},\]
with $\delta \in \{0,1\}$ zero if and only if one of the product $\carmult^{s_i} \carmultbis$ is trivial and with $z$ the number of trivial characters among the $\carmult^{s_j}\carmultbis$.
\end{thm}

\begin{pf}
As the computation is quite lengthly, we do not reproduce it here and refer to \cite[\S 9]{CdlORV.I} and \cite[\S 14 p.~144-150]{CdlORV.II}. The formulas we gave above differ from those given by Candelas, de la Ossa, and Rodriguez Villegas by a factor $\frac{1}{q-1}$ because they compute the number of points in the affine space instead of the projective space. Note also that, although they only showed the formulas for a specific choice of multiplicative and additive $p$-adic characters, it is straightforward to extend their method to arbitrary characters with values in any algebraic closed field of characteristic zero (we only need to replace the relations on the Dwork character by the orthogonality formulas). \qed
\end{pf}

The first thing we do is modify the formulas of Candelas, de la Ossa, and Rodriguez Villegas thanks to the product formula~(\ref{formula:multiplication}).

\begin{lem}
With the notations of Theorem~\ref{result:CdlORV.psi<>0},
\begin{equation}\label{formula:nb.zero.quintic:quotient}
N_{(s_1,s_2,s_3,s_4,s_5),\carmult}(q) = \frac{1}{q-1} \sum_{\carmultbis \in \Fqhat}{q^{3-z-\delta}\frac{\prod_{j=1}^{5}{\gausssum{\carmult^j\carmultbis}}} {\prod_{j=1}^{5}{\gausssum{\carmult^{s_j}\carmultbis}}} \carmultbis(\tfrac{1}{\parameter^5})}.
\end{equation}
\end{lem}

\begin{pf}
In the formula for $N_{(s_1,s_2,s_3,s_4,s_5),\carmult}(q)$ of Theorem~\ref{result:CdlORV.psi<>0}, we write $\carmultbis(\tfrac{1}{(5\parameter)^5}) = \frac{1}{\carmultbis(5)^5} \carmultbis(\tfrac{1}{\parameter^5})$ and replace $\carmultbis(5)^5$ by $q^2\frac{G(\caradd,\carmultbis^5)}{\prod_{j=1}^{5}{G(\caradd,\carmult\carmultbis)}}$ as per formulas~(\ref{formula:multiplication}) and (\ref{formula:product.simplification}). \qed
\end{pf}

The first consequence of this new formula is that $N_\textup{sing}(q)$ is non zero only when $\parameter^5 = 1$ i.e. only when the quintic is singular. In the case $q=p$, this was shown by Candelas, de la Ossa, and Rodriguez Villegas in \cite[\S 9.3]{CdlORV.I} by using the multiplication formula for the $p$-adic Gamma function. We have merely replaced this formula by the multiplication formula~(\ref{formula:multiplication}) for Gauss sums. It both simplifies the proof and allow to treat at the same time the case $q \neq p$.

\begin{lem}\label{result:nb.zeros.sing}
When $q \equiv 1 \mod 5$,
\[N_{(0,1,2,3,4),\carmult}(q) = \begin{cases} 0 & \text{if $\parameter^5 \neq 1$} \\ q^2 & \text{if $\parameter^5 = 1$} \end{cases}\]
and so the term $24 N_\textup{sing}(q)$ does not contribute to the zeta function $\zetafunction{\quintic{\parameter}/\Fq}{t}$ when $\parameter^5 \neq 1$ and gives rise to a factor $\frac{1}{(1-q^2t)^{24}}$ when $\parameter^5 = 1$.
\end{lem}

\begin{pf}
When $(s_1,\dots,s_5) = (0,1,2,3,4)$, we have $z+\delta = 1$ for all values of $\carmultbis$ in formula~(\ref{formula:nb.zero.quintic:quotient}) and the two products simplify completely. Thus,
\[N_{(0,1,2,3,4),\carmult}(q) = \frac{q^{2}}{q-1} \sum_{\carmultbis \in \Fqhat}{\carmultbis(\tfrac{1}{\parameter^5})}.\]
We now use the fact that when $\parameter^5 \neq 1$, $\carmult(\tfrac{1}{\parameter^5})$ is a $\frac{q-1}{5}$-th root of unity $\neq 1$ and so
\[\sum_{\carmultbis \in \Fqhat}{\carmultbis(\tfrac{1}{\parameter^5})} = \begin{cases} 0 & \text{if $\parameter^5 \neq 1$} \\ q-1 & \text{if $\parameter^5 = 1$}\end{cases} \qed\]
\end{pf}

We now proceed to give the link between $N_a(q)$ and $\card{\curveA{\parameter}(\Fq)}$ and between $N_b(q)$ and $\card{\curveB{\parameter}(\Fq)}$.

\begin{lem}\label{result:link.Na.A}
Suppose $q \equiv 1 \mod 5$. With the notations of Corollary~\ref{result:nbzeros.AB:decomp},
\[N_{(0,0,0,1,4),\carmult}(q) = q\nbzeros{\curveA{\parameter}/\Fq,\carmult} \quad \text{ and } \quad N_{(0,0,0,2,3),\carmult}(q) = q\nbzeros{\curveA{\parameter}/\Fq,\carmult^2},\]
and so the term $-10qN_a(q)$ gives rise to a factor $P_{\curveA{\parameter}/\Fq}(qt)^{10}$ in the zeta function $\zetafunction{\quintic{\parameter}/\Fq}{t}$.
\end{lem}

\begin{pf}
Let us for example do the computations for $N_{(0,0,0,1,4),\carmult}(q)$. Using formula~(\ref{formula:nb.zero.quintic:quotient}), we have, doing obvious simplifications,
\[N_{(0,0,0,1,4),\carmult}(q) = \frac{1}{q-1} \sum_{\carmultbis \in \Fqhat}{q^{3-z-\delta}\frac{\gausssum{\carmult^2\carmultbis} \gausssum{\carmult^3\carmultbis}} {\gausssum{\carmultbis}^2} \carmultbis(\tfrac{1}{\parameter^5})}.\]
If, as in Theorem~\ref{result:nbzeros.AB.psi<>0}, $\nu$ is the number of trivial characters in the pair $(\carmultbis,\carmultbis)$ (here, $a=2$ and $b=3$, so $a+b=5$), we easily see that $z+\delta = 1 + \nu$ and so $3 - z - \delta = 2 - \nu$. Thus, using formula~(\ref{result:nbzeros.AB.carmult}) for $\curveA{\parameter}$,
\[N_{(0,0,0,1,4),\carmult}(q) = \frac{1}{q-1} \sum_{\carmultbis \in \Fqhat}{q^{2-\nu}\frac{\gausssum{\carmult^2\carmultbis} \gausssum{\carmult^3\carmultbis}} {\gausssum{\carmultbis}^2} \carmultbis(\tfrac{1}{\parameter^5})} = q\nbzeros{\curveA{\parameter}/\Fq,\carmult}.\]
The proof for $N_{(0,0,0,2,3),\carmult}(q)$ is similar. \qed
\end{pf}

\begin{lem}\label{result:link.Nb.B}
Suppose $q \equiv 1 \mod 5$. With the notations of Corollary~\ref{result:nbzeros.AB:decomp},
\[N_{(0,0,1,1,3),\carmult}(q) = q\nbzeros{\curveB{\parameter}/\Fq,\carmult} \quad \text{ and } \quad N_{(0,0,1,2,2),\carmult}(q) = q\nbzeros{\curveB{\parameter}/\Fq,\carmult^2},\]
and so the term $-15qN_b(q)$ gives rise to a factor $P_{\curveB{\parameter}/\Fq}(qt)^{15}$ in the zeta function $\zetafunction{\quintic{\parameter}/\Fq}{t}$.
\end{lem}

\begin{pf}
The proof is the same than for the previous lemma. \qed
\end{pf}

Lemma~\ref{result:link.Na.A} gives formula~(\ref{aim:Na=NA:psi<>0}) of Theorem~\ref{main.theorem} of the introduction. Lemma~\ref{result:link.Nb.B} gives the analogous formula for $\curveB{\parameter}$. We have thus proved Theorem~\ref{main.theorem} in all cases.

\section{Remarks on the singular case}\label{section:singular.case}

When $\parameter^5 = 1$ and $q \equiv 1 \mod 5$, the only thing that changes by comparison with the non singular case is that $24 N_\textup{sing}(q) = 24q^2$ and so the term $24 N_\textup{sing}(q)$ gives rise to a factor $\frac{1}{(1-q^{2}t)^{24}}$ in $\zetafunction{\quintic{\parameter}/\Fq}{t}$. Moreover, in that case, $N_a(q)$ and $N_b(q)$ have very simple expressions as the equations of $\curveA{\parameter}$ and $\curveB{\parameter}$ become $y^5 = x^2(1-x)^5$ and so, by using the same methods as in \S\ref{section:nb.points.curves:general.remarks} and \S\ref{section:nb.points.curves:parameter=0} for the case $\parameter = 0$,
\[\card{\curveA{\parameter}(\Fq)} = \card{\curveB{\parameter}(\Fq)} = \begin{cases} q-4 & \text{if $q \equiv 1 \mod 5$} \\ q & \text{otherwise} \end{cases}\]
The zeta function of $\curveA{\parameter}$ and $\curveB{\parameter}$ over $\Fq$ is thus $\frac{(1-t)^{4}}{1-qt}$ when $q \equiv 1 \mod 5$. This result was also shown in \cite[\S 11.3 p.~131-132]{CdlORV.II} by using a change of variable.

For the sake of completeness, let us mention that when $\parameter^5 = 1$, the factor $P_w$ has a modular origin. More precisely, using results of Schoen \cite[p.~109-110]{Schoen.quintic}, Candelas, de la Ossa, and Rodriguez Villegas \cite[\S 12 p.~133-134]{CdlORV.II} showed that, when $\parameter^5 = 1$ and $q=p$,
\[P_w(t) = (1-(\tfrac{5}{p}) pt)(1-a_pt+p^3t^2),\]
where $(\tfrac{\cdot}{p})$ is the quadratic character of $\Fp^*$ (Legendre symbol) and $a_p$ the $p$-th Fourier coefficient of a weight $4$ and level $25$ modular form that can be computed explicitly thanks to the Dedekind $\eta$ function (see \cite[\S 12 p.~134]{CdlORV.II})
\[\begin{split}
f(q)
& = \eta(q^5)^4 \Big( \eta(q)^4 + 5 \eta(q)^3 \eta(q^{25}) + 20 \eta(q)^2 \eta(q^{25})^2 + 25 \eta(q)\eta(q^{25})^3 + 25 \eta(q^{25})^4 \Big)\\
& = q + q^2 + 7q^3 - 7q^4 + 7q^6 + 6q^7 - 15q^8 + 22q^9 - 43q^{11} -
49 q^{12} - 28 q^{13} \\  
& \quad + 6 q^{14} + 41 q^{16} + 91 q^{17} + 22 q^{18} - 35 q^{19} + 42 q^{21} - 43 q^{22} + 162 q^{23} - 105 q^{24} \\  
& \quad - 28 q^{26} - 35 q^{27} - 42 q^{28} + 160 q^{29} + 42 q^{31} + 161 q^{32} - 301 q^{33} + 91 q^{34} \\  
& \quad - 154 q^{36} - 314 q^{37} - 35 q^{38} - 196 q^{39} - 203 q^{41} + 42 q^{42} + 92 q^{43} + 301 q^{44} \\  
& \quad + 162 q^{46} + 196 q^{47} + 287 q^{48} - 307 q^{49} + 637 q^{51} + 196 q^{52} + 82 q^{53} + \dots 
\end{split}\]

\begin{ack}
The author would like to thank his advisor, J.~Oesterl\'e, for his comments on preliminary versions of this paper.
\end{ack}


\begin{thebibliography}{00}

\bibitem{BEW}
B.~C. Berndt, R.~J. Evans, K.~S. Williams, Gauss and Jacobi Sums, vol.~21 of Canadian Mathematical Society Series of Monographs and Advanced Texts, Willey Interscience, 1998.

\bibitem{CdlORV.I}
P.~Candelas, X.~de~la Ossa, F.~Rodriguez~Villegas, Calabi Yau manifolds over finite fields, {I}, preprint. Available online at \url{http://www.arxiv.org/hep-th/0012233}.

\bibitem{CdlORV.II}
P.~Candelas, X.~de~la Ossa, F.~Rodriguez~Villegas, Calabi Yau manifolds over finite fields, {II}, in: N.~Yui, J.~D. Lewis (eds.), Proceedings of the Workshop on ``Calabi {Y}au Varieties and Mirror Symmetry'', Fields Institute, Toronto, July 23--29, 2001, vol.~38 of Fields Inst. Comm. Series, Amer. Math. Soc., 2003, 121--157.

\bibitem{GY}
F.~Gouvea, N.~Yui, Arithmetic of Diagonal Hypersurfaces over Finite Fields, London Mathematical Society Lecture Note Series 209, Cambridge University Press.

\bibitem{GR}
B.~Gross, D.~Rohrlich, Some Results on the Mordell-Weil Group of the Jacobian of the Fermat Curve, Invent. Math. 44 (1978) 201--224.

\bibitem{KOB.nb.points.hyper}
N.~Koblitz, The number of points on certain families of hypersurfaces over finite fields, Compositio Mathematica 48 (1983) 3--23.

\bibitem{Schoen.quintic}
C.~Schoen, On the geometry of a special determinantal hypersurface associated to the Mumford-Horrocks vector bundle, J. Reine Angew. Math. 364 (1986) 85--111.

\bibitem{WEIL.nb.sol.finite.fields}
A.~Weil, Numbers of solutions of equations in finite fields, Bull. Amer. Math. Soc. 55 (1949) 497--508.

\end{thebibliography}







\end{document}